\documentclass[review,12pt]{elsarticle}

\usepackage{amssymb}

\journal{Communications in Nonlinear Science and Numerical Simulation}

\usepackage{xspace}
\usepackage{algorithm}
\usepackage{algorithmic}
\usepackage{subfigure}
\usepackage{amsfonts}
\usepackage{amsmath}
\usepackage{enumerate}
\usepackage{framed}
\usepackage{graphicx}
\usepackage{lscape}
\usepackage{bm}
\usepackage{color}
\usepackage{multicol}
\usepackage{amsthm}

\newcommand{\figref}[1]{{Figure~\ref{#1}}}

\newcommand{\secref}[1]{{Section~\ref{#1}}}
\renewcommand{\eqref}[1]{{(\ref{#1})}}

\begin{document}

\begin{frontmatter}


\title{Localized numerical impulse solutions in diffuse neural networks modeled by the complex fractional Ginzburg-Landau equation}

\author[yd1,am,cetic]{Alain Mvogo\corref{cor1}}
\ead{mvogo@aims.ac.za} \cortext[cor1]{Corresponding author}
\address[yd1]{Department of Physics, Faculty of Science, University of Yaounde I, P.O. Box 812, University of Yaounde, Cameroon}
\address[am]{The African Institute for Mathematical Sciences (AIMS),
6-8 Melrose Rd, Muizenberg 7945, South Africa}
\address[cetic]{Centre d'Excellence Africain en Technologies de l'Information et de la Communication, University of Yaounde I, Cameroon}
\author[at,atb]{Antoine Tambue}
\ead{antonio@aims.ac.za}
\address[at]{The African Institute for Mathematical Sciences (AIMS) and Stellenbosch University, 6-8 Melrose Road, Muizenberg 7945, South Africa}
\address[atb]{Center for Research in Computational and Applied Mechanics (CERECAM), and Department of Mathematics and Applied Mathematics, University of Cape Town, 7701 Rondebosch, South Africa.}
\author[yd1,cetic]{Germain H. Ben-Bolie}
\ead{gbenbolie@yahoo.fr}
\author[yd1,cetic]{Timol\'{e}on C. Kofan\'{e}}
\ead{tckofane@yahoo.com}
\address{}
\begin{abstract}

We investigate localized wave solutions in a network of
Hindmarsh-Rose neural model taking into account the long-range
diffusive couplings. We show by a specific analytical technique that
the model equations in the infrared limit (wave number $k\rightarrow
0$) can be governed by the complex fractional Ginzburg-Landau (CFGL)
equation. According to the stiffness of the system, we propose both
the semi and the linearly implicit Riesz fractional
finite-difference schemes to solve efficiently the CFGL equation.
The obtained fractional numerical solutions for the nerve impulse
reveal localized short impulse properties. We also show  the
equivalence between the continuous CFGL and the  discrete
Hindmarsh-Rose models for relatively large network.\\
\end{abstract}

\end{frontmatter}

\section{Introduction}
\label{sec:intr}

The interest in investigating and controlling the propagation of
waves in neural tissues has been increasingly growing during the
last decades. This is because the conditions under which cortical
waves occur are very primordial in the understanding of the normal
processing of sensory stimuli as well as more pathological forms of
behavior \cite{Kilpatrick,Clarke}. In that sense,  many studies have
been carried out that indicated the presence of localized nonlinear
waves in the neural systems (see, e.g., \cite{Kakmeni} and
references therein). Interestingly, the recent work by Kakmeni
\emph{et al.} reported on the presence of these waves of the nerve
impulse in diffusive Hindmarsh-Rose (HR) neural networks with
nearest-neighbor couplings \cite{Kakmeni}. Also, as recently
demonstrated, the dynamics of an individual neuron in diffusive HR
neural networks may be influenced by the interaction or coupling
with other neurons \cite{Djeundam,Steur}. It would be interesting
now to see what are the effects that the long-range diffusive
coupling or interaction has on the wave propagation in such a
network.

The intuitively obvious fact that many biological systems are
systems with memory or aftereffects is now confirmed by many
researches. The modeling of these systems by fractional-order
differential equations has more advantages than the classical
mathematical modeling using the integer-order, in which such effects
are neglected. As it has been shown, even processing of external
stimuli by individual neural oscillator can be described by
fractional differentiation \cite{Lund1,Lund2}. In many cases memory
effect obeys  the power law and the corresponding system could be
described by fractional differential equation. It becomes also now
interesting  to investigate localized waves in such an equation in
diffusive neural networks when the long-range coupling is taken into
account.

We aim in this paper to study the properties of localized waves in
the diffusive HR neural networks with long-range interactions that
can work in some way as a long memory. We  show  that the model can
be governed in the infrared limit (wave number $k\rightarrow 0$) can
be governed by the complex fractional Ginzburg-Landau (CFGL)
equation. According to the stiffness of the system, we propose both
the semi and the linearly implicit Riesz fractional
finite-difference schemes to solve efficiently the CFGL equation.
The obtained fractional numerical solutions for the nerve impulse
reveal localized short impulse properties. The fractional order
mostly contributes to the behavior of the tails of the impulse. It
is also shown  the equivalence between the continuous CFGL and the
discrete Hindmarsh-Rose models for relatively large network.

The rest of the paper is organized as follows. In \secref{sec2}, we
present the neural network taking into account the long-range
diffusive coupling. In \secref{sec3}, by means of the perturbation
technique, we derive the complex fractional Ginzburg-Landau (CFGL)
equation which describes the equation of motion. In \secref{sec4},
we solve efficiently the CFGL equation following closely \cite{shen}
in space discretization, and propose the semi-implicit Riesz
fractional finite-difference scheme and the linearly implicit Riesz
fractional finite-difference scheme where only one linear system is
solved by time iteration. We then present the numerical results and
show  the equivalence between the continuous CFGL model and the
discrete  HR model for relatively large network. Our work is
summarized in \secref{sec5}.


\section{The Hindmarsh-Rose coupled model}
\label{sec2}

Many nontrivial examples of dynamical systems have been provided by
phenomenological and neurophysiological models developed to
reproduce the activities of neural oscillators. The Hindmarsh-Rose
model \cite{Hindmarsh}, a generalization of the Fitzhugh equations
\cite{Fitzhugh}, represents a paradigmatic example of these systems.
It aims to study the spiking-bursting behavior of the membrane
potential observed in the single neuron experiments. In this paper,
following Refs. \cite{Kakmeni,Djeundam} we generalize the HR neural
model assuming only that the coupling between neural oscillators are
long-ranged through the membrane potential variable. The HR neural
network is then assumed as a system of $N$ neural oscillators in
which the configuration of couplings is assumed to be power
long-ranged. In this case, each unit of HR neural model is coupled
to any other. The model can be reformulated by means of the
following nonlinear ordinary differential equations:

\begin{eqnarray}
\label{1} \left\lbrace \begin{array}{l}
\dot{u}_{n}=v_{n}-au_{n}^{3}+bu_{n}^{2}-w_{n}+I
+\sum\limits_{m=1,m\neq n}^{N} K_{\alpha}(n-m)(u_{n}-u_{m}),\\
\dot{v}_{n}=c-du_{n}^{2}-ev_{n},\\
\dot{w}_{n}=r[s(u_{n}-u_{0})-w_{n}],
\end{array}\right.
\end{eqnarray}
where the variable $u$ is the membrane potential (nerve impulse),
$v$ is the spiking variable which takes into account the measure of
the rate at which transport of sodium and potassium ions is made
through fast ion channels, and $w$ is the bursting variable which
takes into account the rate at which the transport of other ions
($Cl^{-}$ and proteins anions) made through slow ions channels. The
values of the parameters of the HR model are $a = 1.0$, $b = 3.0$,
$c = 1.0$, $d = 5.0$, $r = 0.008$, $s = 4.0$, $e = 1.0$, $u_{0} =
-1.60$ and $2.92 < I < 3.40$.

Physiologically, responses generated within the cell can travel not
only to neighboring cells through intercellular communication using
a gap-junction but also through extracellular communication,
involving the secretion of molecular signals such as
neurotransmitters. Then, in comparison to the model of Ref.
\cite{Kakmeni}, the present HR model is generalized via the presence
of the term $\sum\limits_{m\neq n} K_{\alpha}(m-n)(u_{n}-u_{m})$
which characterizes long-range diffusive interaction in the system
and can appear as an effective interaction in dispersive and complex
systems \cite{Zosimov,Mvogo1,Mvogo2}. The latter is due to the fact
that the extracellular messenger can propagate from one cell to its
direct neighbors  and even extend to other neighboring
non-contacting cells. The nonlocal coupling interaction is given by
the power-law dependence
\begin{equation}
\label{2} K_{\alpha}(n)=\frac{K}{|n|^{\alpha+1}},
\end{equation}
where $K$ is the coupling parameter such as the synaptic strength,
while $\alpha$ which is the LRI parameter, physically describes a
level of collective interaction of neural oscillators.

It has been demonstrated by Steur \emph{et al.}  \cite{Steur} that
such a neural system (1) coupled via diffusive coupling is
semi-passive, then the solutions of all connected systems in the
network are bounded. We are interested by nonlinear waves in the
network.  To achieve this, we first differentiate the first equation
of \eqref{1} and substitute $\dot{v}_{n}$ into the obtained
second-order ordinary differential equation. Then, we rewrite
suitably \eqref{1} in a Lienard form, that is a second-order
differential equation with a small damping term, such that

\begin{eqnarray}
\label{3} \left\lbrace \begin{array}{l}
\ddot{u}_{n}+\Omega_{0}^{2}u_{n}+(\eta_{0}+\eta_{1}u_{n}+\eta_{2}u_{n}^{2})\dot{u_{n}}+\lambda_{1}u_{n}^{2}+\frac{\eta_{2}}{3}u_{n}^{3}+\lambda_{3}w_{n}+I_{0}\\
\newline
=\sum\limits_{m=1,m\neq n}^{N}\frac{1}{|n-m|^{\alpha+1}}\Big[
c_{0}(u_{n}-u_{m})+c_{1}(\dot{u}_{n}-\dot{u}_{m})\Big]\\
\dot{v_{n}}=c-du_{n}^{2}-ev_{n}\\
\dot{w_{n}}=r[s(u_{n}-u_{0})-w_{n}],
\end{array}\right.
\end{eqnarray}
where $\Omega_{0}$, $\eta_{0}$, $\eta_{1}$, $\eta_{2}$,
$\lambda_{1}$, $\lambda_{3}$, $I_{0}$, $c_{0}$ and $c_{1}$ are
constant parameters related to those of \eqref{1} such as
$\Omega_{0}^{2}=rs$, $\eta_{0}=e$, $\eta_{1}=2b$, $\eta_{2}=3ea$,
$\lambda_{1}=d-eb$, $\lambda_{3}=e-r$, $I_{0}=c+rsu_{0}+I$,
$c_{0}=eK$ and $c_{1}=K$. In general, the solutions of \eqref{3} can
be obtained using perturbation techniques. In that sense, we
introduce the following variables $u_{n}=\varepsilon \psi_{n}$,
$v_{n}=\varepsilon \Phi_{n}$ and $w_{n}=\varepsilon \beta_{n}$,
where $\varepsilon<<1$. By keeping in the development the first two
nonlinear terms, the governing equations of motion in the neural
network then become

\begin{eqnarray}
\label{4} \left\lbrace \begin{array}{l}
\ddot{\psi}_{n}+\Omega_{0}^{2}\psi_{n}+\varepsilon(\varepsilon\eta_{0}+\eta_{1}\psi_{n}+\varepsilon
\eta_{2}\psi_{n}^{2})\dot{\psi}_{n}+\varepsilon\lambda_{1}\psi_{n}^{2}+\varepsilon^{2}\frac{\eta_{2}}{3}\psi_{n}^{3}+\varepsilon^{2}\lambda_{3}\beta_{n}\\
\newline
=\sum\limits_{m=1,m\neq n}^{N}\frac{1}{|n-m|^{\alpha+1}}\Big[
c_{0}(\psi_{n}-\psi_{m})+\varepsilon^{2}c_{1}(\dot{\psi}_{n}-\dot{\psi}_{m})\Big]\\
\dot{\Phi}_{n}+\varepsilon d\psi_{n}^{2}+e\Phi_{n}=0\\
 \dot{\beta}_{n}+r\beta_{n}-\Omega_{0}^{2}\psi_{n}=0.
\end{array}\right.
\end{eqnarray}

While writing \eqref{4}, the coupling parameter $\lambda_{3}$ of the
membrane potential with the bursting variable has been perturbed of
order $\varepsilon^{2}$, taking into account the fact that the
variation of the bursting variable is slower than the one of the
membrane potential. In addition, as we are interested by an analysis
in a weakly dissipative medium, we have assumed the parameters
$\eta_{0}$ and $c_{1}$ to be perturbed at the order
$\varepsilon^{2}$.

\section{Equation of motion} \label{sec3}

The chain of the neural network is a very long having several
thousand neural oscillators compared to the distance between the
neighboring neural oscillators along the chain. It is therefore
appropriate to make a continuum approximation, which is also valid
in the long wavelength limit. The non-locality features of the
medium often impose the necessity of using non-traditional tools. In
that follows, we first assume for Eq. \eqref{4} the following
solutions

\begin{eqnarray}
\label{5} \left\lbrace \begin{array}{l}
 \psi_{n}=\varepsilon(B_{n}^{(1)}e^{\mathrm{i}\theta_{n}}+c.c)+\varepsilon^{2}[C_{n}^{(1)}+(D_{n}^{(1)}e^{\mathrm{2i}\theta_{n}}+c.c)]\\
\label{5}
\Phi_{n}=\varepsilon(B_{n}^{(2)}e^{\mathrm{i}\theta_{n}}+c.c)+\varepsilon^{2}[C_{n}^{(2)}+(D_{n}^{(2)}e^{\mathrm{2i}\theta_{n}}+
c.c)]\\
\label{5}
\beta_{n}=\varepsilon(B_{n}^{(3)}e^{\mathrm{i}\theta_{n}}+c.c)+\varepsilon^{2}[C_{n}^{(3)}+(D_{n}^{(3)}e^{\mathrm{2i}\theta_{n}}+
c.c)],
\end{array}\right.
\end{eqnarray}
with $\theta_{n}=kn-\Omega t$, where $k$ is the normal mode wave
vector and $\Omega$ is the angular velocity of the wave. The
variable $t$ is rescaled through the perturbative small parameter
$\varepsilon$ as $t\rightarrow\varepsilon^{2}t$.

In the following, we replace the  solutions (\ref{5}) and their
derivatives in the
 new membrane potential equation of motion
given by the first equation of \eqref{4}. We then group the terms in
the same power of $\varepsilon$, which leads us to a system of
equations. Each of these equations will correspond to each
approximation for specific harmonics. To reach this goal, we
consider the infinite network of neural oscillators ($N\rightarrow
\infty$). We multiply Eq. (\ref{4}) by $\exp(-ikn)$ and we sum over
$n$ from $-\infty$ to $+\infty$. Then, we introduce the following
functions

\begin{eqnarray}
\label{6} \left\lbrace \begin{array}{l} f^{j}(k,
t)=\sum\limits_{n=-\infty}^{\infty} e^{-\mathrm{i}kn}B_{n}^{(j)}(t)\\
 g^{j}(k,t)= \sum\limits_{n=-\infty}^{\infty} e^{-\mathrm{i}kn}C_{n}^{(j)}(t)\\
h^{j}(k,t)= \sum\limits_{n=-\infty}^{\infty} e^{-\mathrm{i}kn}D_{n}^{(j)}(t)\\
\tilde{J}_{\alpha}(k)=\sum\limits_{n=-\infty}^{\infty}e^{-\mathrm{i}k
n}\frac{1}{|n|^{\alpha+1}},
\end{array}\right.
\end{eqnarray}
with $j=1,2,3$,  and

$$\tilde{J}_{\alpha}(0)=2\sum_{n=1}^{\infty}\frac{1}{|n|^{\alpha+1}}=2\zeta(\alpha+1),$$
where $\zeta(\alpha)$ is the Riemann zeta function.

In the long-wave limit, we may adopt $f^{j}(k, t)$, $g^{j}(k, t)$
and  $h^{j}(k, t)$  as $k^{th}$ Fourier components of continuous
functions $B^{(j)}(x, t)$, $C^{(j)}(x, t)$ and $D^{(j)}(x, t)$,
respectively such that $B_{n}^{(j)}(t)\rightarrow B^{(j)}(x, t)$,
$C_{n}^{(j)}(t)\rightarrow C^{(j)}(x, t)$ and
$D_{n}^{(j)}(t)\rightarrow D^{(j)}(x, t)$. The functions are related
each other by the Fourier transforms such that

\begin{eqnarray}
\label{FT} \left\lbrace \begin{array}{l}
B^{(j)}(x,t)=\frac{1}{2\pi}\int_{-\infty}^{\infty}e^{\mathrm{i}kx}f^{j}(k,t)dk\\
 C^{(j)}(x,t)=\frac{1}{2\pi}\int_{-\infty}^{\infty}e^{\mathrm{i}kx}g^{j}(k,t)dk\\
D^{(j)}(x,t)=\frac{1}{2\pi}\int_{-\infty}^{\infty}e^{\mathrm{i}kx}h^{j}(k,t)dk.
\end{array}\right.
\end{eqnarray}

\begin{figure}[h!]
\centering
\includegraphics[width=3in]{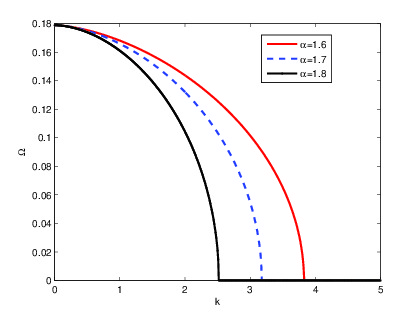}
\caption{The dispersion relation of the nerve impulse according to
three different values of $\alpha$, namely $\alpha=1.6$,
$\alpha=1.7$ and $\alpha=1.8$. $\Omega_{0}^{2}=0.032$ and
$c_{0}=0.001$. It appears that the long-rang parameter $\alpha$
affects the dispersion area.}\label{fig1}
\end{figure}

After some algebras, at the order $\varepsilon^{1}$ after the
annihilation of terms in $e^{\pm \mathrm{i}\theta}$, we obtain the
relation

\begin{equation}\label{9}
\Omega^{2}=\Omega_{0}^{2}+c_{0}a_{\alpha}|k|^{\alpha},
\end{equation}
which determines the dispersion relation of linear waves of the
system. As displayed in \figref{fig1}, the corresponding linear spectrum is
reduced when $\alpha$ increases. To obtain \eqref{9}, we have used
the infrared approximation \cite{Tarasov,Laskin}

\begin{equation}\label{10}
[\tilde{J}_{\alpha}(0)-\tilde{J}_{\alpha}(k)]\approx
a_{\alpha}|k|^{\alpha},\qquad (0<\alpha<2, \alpha\neq1)
\end{equation}
where $a_{\alpha}=2\Gamma(-\alpha)\cos(\pi\alpha/2)$.

At the order $\varepsilon^{2}$, terms without exponential dependence
give

\begin{equation}\label{8}
C^{(1)}=-\frac{2\lambda_{1}}{\Omega_{0}^{2}}|B^{(1)}|^{2},
\end{equation}
while at the same order, terms with $e^{2\mathrm{i}\theta}$ give the
relation

\begin{equation}\label{11}
D^{(1)}=\frac{\lambda_{1}-\mathrm{i}\Omega\eta_{1}}{\Omega_{1}^{2}}(B^{(1)})^{2}.
\end{equation}
where
$\Omega_{1}^{2}=2\Omega_{0}^{2}+\Omega^{2}+4(\Omega^{2}-\Omega_{0}^{2})^{2}$.\\

For the third equation of \eqref{3}, at the order $\varepsilon^{1}$
the terms with $e^{\mathrm{i}\theta}$ give the relation
\begin{equation}\label{12}
B^{(3)}=\frac{\Omega_{0}^{2}(r+\mathrm{i}\Omega)B^{(1)}}{r^{2}+\Omega^{2}}.
\end{equation}

Collecting all the terms depending on $e^{\mathrm{i}\theta}$ in
\eqref{4} at the order $\varepsilon^{3}$, we obtain the following
equation

\begin{equation}\label{13}
\begin{split}
-2\mathrm{i}\Omega\frac{\partial B^{(1)}}{\partial
t}= &\mathrm{i}\Omega\eta_{0}B^{(1)}+(\mathrm{i}\Omega\eta_{1}-2\lambda_{1})(B^{(1)}C^{(1)}+B^{(1*)}D^{(1)})\\
&+(\mathrm{i}\Omega-1)\eta_{2}|B^{(1)}|^{2}B^{(1)}-\lambda_{3}\frac{\Omega_{0}^{2}(r+\mathrm{i}\Omega)}{r^{2}+\Omega^{2}}B^{(1)}\\
&-\mathrm{i}c_{1}\Omega (a_{\alpha}|k|^{\alpha}B^{(1)}).
\end{split}
\end{equation}

Rewriting this equation taking into account  the connection between
the  Riesz fractional derivative  and its Fourier transform
\cite{Samko}

\begin{equation}\label{14}
|k|^{\alpha}\longleftrightarrow-\frac{\partial^{\alpha}}{\partial
|x|^{\alpha}}, \;\
|k|^{2}\longleftrightarrow-\frac{\partial^{2}}{\partial |x|^{2}}
\end{equation}
we obtain

\begin{equation}\label{15}
\frac{\partial B^{(1)}}{\partial t}=\gamma
B^{(1)}+P_{r}\frac{\partial^{\alpha}B^{(1)}}{\partial
|x|^{\alpha}}-Q|B^{(1)}|^{2}B^{(1)},
\end{equation}
where the coefficients $\gamma$, $P_{r}$ and $Q$
are given by\\
$\gamma=\gamma_{r}+\mathrm{i}\gamma_{i}$, $P_{r}=c_{1}a_{\alpha}/2$
and  $Q=Q_{r}+\mathrm{i}Q_{i}$.

The coefficients $\gamma_{r}$ and $\gamma_{i}$ are the real and
imaginary parts of the dissipation coefficient. For the nonlinearity
coefficient the same terminology is used. The coefficients $Q_{r}$,
$Q_{i}$, $\gamma_{r}$ and $\gamma_{i}$ are given by
\begin{eqnarray*}
\begin{split}
&\gamma_{r}=\dfrac{\lambda_{3}\Omega_{0}^{2}}{2(r^{2}+\Omega^{2})}-\dfrac{\eta_{0}}{2},\,\,\,\,\,
\gamma_{i}=-\dfrac{r\lambda_{3}\Omega_{0}^{2}}{2\Omega(r^{2}+\Omega^{2})},\\
&Q_{r}=\dfrac{\eta_{2}}{2}-\dfrac{\eta_{1}\lambda_{1}}{\Omega_{0}^{2}}+\dfrac{3\eta_{1}\lambda_{1}}{2\Omega_{1}^{2}},\,\,\,\,\,
 Q_{i}=\dfrac{1}{\Omega}
\left(\dfrac{\eta_{2}}{2}-\dfrac{\Omega^{2}\eta_{1}^{2}-2\lambda_{1}^{2}}{2\Omega_{1}^{2}}-\dfrac{2\lambda_{1}^{2}}{\Omega_{0}^{2}}\right).
\end{split}
\end{eqnarray*}

Equation (\ref{15}) which is a new general theoretical framework
derived in our neural network is the complex fractional
Ginzburg-Landau equation. This confirm the fact that the brain may
actively work effectively using the spatial dimension for
information processing but not only in time domain
\cite{Negro,Pedro,Azad}. In Ref. \cite{Tarasov2}, the fractional
Ginzburg-Landau equation is derived from the variational Euler
Lagrange equation for fractal media. In the present work, we confirm
once more the fact that using the Fourier transforms and the
infrared limit, the long-range interactions lead under special
conditions to the fractional dynamics \cite{Tarasov,Laskin,Samko}.
The fractional Ginzburg-Landau equation has been proposed by
Weitzner and Zaslavsky \cite{Weitzner} to describe the dynamical
processes in a medium with fractal dispersion. Its generalization
has been used by Milovanov and Rasmussen \cite{Milovanov} as an
unconventional approach to critical phenomena in complex media.

In this work, the motion of modulated waves in diffuse neural
networks are proven to be described by the CFGL equation. The
infrared limit of an infinite chain of neural oscillators with the
long-range diffusive interactions can be described by equations with
the fractional Riesz coordinate derivative of order $\alpha<2$. To
the best of our knowledge, this is the first research work that
attempts to describe the dynamical behavior of neural networks with
an equation of fractional order. This result suggests that neurons
can participate in a collective processing of long-scale
information, a relevant part of which is shared over all neurons.

\section{The Semi-implicit Riesz fractional
finite-difference and the linearly Riesz fractional
finite-difference schemes}
 \label{sec4}

In the previous section, we have demonstrated that the HR neural
network can be elegantly described by the CFGL equation. In general,
analytical and closed solutions of fractional equations cannot  be
obtained. In that case, numerical techniques are used to identify
the solution behavior of such fractional equations. In this section,
we provide the semi-implicit Riesz fractional finite-difference
scheme and the linearly implicit Riesz fractional finite-difference
scheme  to find numerically localized wave solutions for the CFGL
equation (\ref{11}). We also show numerically  the equivalent
between the continuous  CFGL model and the discrete  HR  model for
relatively large  network.

To begin our numerical analysis, it is convenient to recall  the
CFGL equation (\ref{15}) as
\begin{equation}
\label{new}
\left\lbrace \begin{array}{l}
\dfrac{\partial B^{1}}{\partial
t}=P_{r}\dfrac{\partial^{\alpha}B^{1}}{\partial
|x|^{\alpha}} + (\gamma_{r}+\mathrm{i}\gamma_{i})
B^{1}-(Q_{r}+\mathrm{i}Q_{i}) |B^{1}|^{2}B^{1}, \\
B^{1}(x,0)=g(x),\,\,\,\,  x \in [0,b_{0}],\\
B^{1}(0,t)= \phi_1(t),\, B^{1}(b_{0},t)=\phi_2(t), \;\,\;\; t\in (0,T],
\end{array}\right.
\end{equation}
The functions $g, \phi_1$ and $ \phi_2$ are sufficiently smooth functions.
Note that $B^{1}\equiv B^{(1)}$, the function
$g$ is the initial solution, $T>0$ is  the final time and
$\dfrac{\partial^{\alpha} }{\partial |x|^{\alpha}}$ is space Riesz
fractional derivative of order $\alpha$ given  for $ 0< \alpha <
2,\, \alpha\neq 1$ by
\begin{eqnarray}
 \dfrac{\partial^{\alpha} }{\partial |x|^{\alpha}}
 = -c_{\alpha} \left(_{-\infty}D_{x}^{\alpha}+ _{x}D_{+ \infty}^{\alpha}\right),
\end{eqnarray}
where  the coefficient
\begin{eqnarray}
\label{new1}
\left\lbrace \begin{array}{l}
c_{\alpha}= \dfrac{1}{2 \cos(\alpha \pi/2)},\\
\newline
_{-\infty}D_{x}^{\alpha}B^{1}= \left(\dfrac{d}{dx}\right)^{m} \left[ _{\infty}I_{x}^{m-\alpha} B^{1}(x,t) \right],\\
_{x}D_{+ \infty}^{\alpha}B^{1}= (-1)^{m} \left(\dfrac{d}{dx}\right)^{m} \left[ _{x}I_{+ \infty}^{m-\alpha} B^{1}(x,t) \right],\\
\end{array}\right.
\end{eqnarray}
with $m \in \mathbb{N}$ such that $m-1 < \alpha \leq m$.  The
terms $_{-\infty}D_{x}^{\alpha}$ and $_{x}D_{+ \infty}^{\alpha}$
are respectively the left and  the right side  Riemann-Liouville
fractional derivatives. The  left and  right side Weyl fractional
integrals used in  \eqref{new1}  are defined by
\begin{eqnarray}
 \left\lbrace \begin{array}{l}
 _{\infty}I_{x}^{m-\alpha} B^{1}(x,t)=\dfrac{1}{ \Gamma(\alpha)} \int_{-\infty}^{x} \left(x-\zeta \right)^{\alpha-1}  B^{1}(\zeta,t) d \zeta,\\
 _{x}I_{+ \infty}^{m-\alpha} B^{1}(x,t)=\dfrac{1}{ \Gamma(\alpha)} \int_{x}^{+ \infty} \left(x-\zeta \right)^{\alpha-1}  B^{1}(\zeta,t) d \zeta.\\
       \end{array}\right.
\end{eqnarray}
By setting $B^{1}=U+ \mathrm{i}V$, where $U$ and $V$ are respectively the real and imaginary parts of $B^{1}$, \eqref{new} is equivalent to the following coupled system

\begin{eqnarray}
\label{new11}
\left\lbrace \begin{array}{l}
\dfrac{\partial U}{\partial
t}=   P_{r}\dfrac{\partial^{\alpha}U}{\partial
\vert x \vert^{ \alpha}} + \gamma_{r} U- \gamma_{i}V - \left( Q_{r} U - Q_{i} V \right) \left( U^{2}+ V^{2}\right),\\
\newline
\dfrac{\partial V}{\partial
t}=   P_{r}\dfrac{\partial^{\alpha}V}{\partial
 \vert x \vert^{\alpha}} + \gamma_{i} U +\gamma_{r}V - \left( Q_{i} U + Q_{r} V \right) \left( U^{2}+ V^{2}\right),\\
\newline
 (U(x,0), V(x,0)= g (x), \;\,\,\,\,\,\, x \in [0,b_{0}],\\
 \newline
U(0,t)= \text{Re}(\phi_1(t)),\, U(b_{0},t)= \text{Re}(\phi_2(t)),\\
 V(0,t)= \text{Im}(\phi_1(t))\;\,\; V(b_{0},t)= \text{Im}(\phi_2(t)), \;\,\;\; t\in (0,T],\\
\end{array}\right.
\end{eqnarray}
where  $\text{Re}$ and $ \text{Im}$ are respectively the real part and  the imaginary part.
Let us use  the following identification $B^{1}\equiv (U, V)^{T}$.
By setting
\begin{eqnarray}
\label{new112}
  \mathcal{A}= \left(\begin{smallmatrix}  P_{r}\dfrac{\partial^{\alpha} }{\partial
\vert x \vert^{\alpha}}+\gamma_{r} & -\gamma_{i}\\ \gamma_{i} & P_{r}\dfrac{\partial^{\alpha} }{\partial
\vert x \vert^{\alpha}}+\gamma_{r} \end{smallmatrix} \right), \\
 F_1(B^{1})= \left( \begin{smallmatrix} - \left( Q_{r} U - Q_{i} V \right) \left( U^{2}+ V^{2}\right)\\ - \left( Q_{i} U + Q_{r} V \right) \left( U^{2}+ V^{2}\right) \end{smallmatrix}\right),
\end{eqnarray}
the coupled system \eqref{new11} becomes
\begin{eqnarray}
 \label{compac}
 \left\lbrace \begin{array}{l}
 \dfrac{\partial B^{1}}{\partial t}=\mathcal{A} B^{1}+ F_1(B^{1}),\\
 B^{1}(x,0)=g(x), \;\,\,\,\,\,\, x \in [0,b_{0}],\\
 B^{1}(0,t)= \phi_1(t), B^{1}(b_{0},t)=\phi_2(t).
 \end{array}\right.
\end{eqnarray}
For space discretization, we use  the weighted Riesz fractional
finite-difference approximation as presented in \cite{shen,inhomo}. We
divide the interval $(0,b_{0})$  into  $M$ sub-interval with the
step $h=b_{0}/M$. In order to perform the space discretization
with our homogeneous boundary conditions, the function $ B^{1}$
should be extended   to the whole $\mathbb{R}$ (see \cite{shen})
as
\begin{eqnarray}
 \label{extended}
B_{*}^{1}(x,t) \equiv (U^{*}, V^{*})^{T}=\left\lbrace \begin{array}{l}
                            B^{1}(x,t),\,\,\,\,\,\, x \in [0,b_0],\\
                           (0,0)^T, \,\,\,\,\,\, x \in (-\infty,0) \cup  (b_0,+\infty).\\
               \end{array} \right.
\end{eqnarray}
Using \eqref{extended}, each component of the function $ B^{1}$ can be discretized  by the centered finite difference as follows for $ 0< \alpha < 2,\, \alpha \neq 1$
\begin{eqnarray}
\label{dis1}
\dfrac{\partial^{\alpha}  B_{*}^{1} }{\partial|x|^{\alpha}}=
\left(\begin{smallmatrix}
 -\dfrac{1}{h^{\alpha}} \sum\limits_{- \infty}^{+ \infty}
w_{k}^{\alpha} U^{*}(x-kh,t) + \mathcal{O}(h^{2})\\
-\dfrac{1}{h^{\alpha}}\sum\limits_{- \infty}^{+ \infty}
w_{k}^{\alpha} V^{*}(x-kh,t) + \mathcal{O}(h^{2})
\end{smallmatrix} \right).
\end{eqnarray}
Since $ B_{*}^{1}(x,t)=(0,0)^{T}$  for $  x \in (-\infty,0) \cup  (b_0,+\infty)$,  we therefore have
\begin{eqnarray}
\label{dis1}
\dfrac{\partial^{\alpha}  B^{1} }{\partial|x|^{\alpha}}=
\left(\begin{smallmatrix}
 -\dfrac{1}{h^{\alpha}} \sum\limits_{-(b_{0}-x)/h}^{(x-0)/h}
w_{k}^{\alpha} U(x-kh,t) + \mathcal{O}(h^{2})\\
-\dfrac{1}{h^{\alpha}}\sum\limits_{-(b_{0}-x)/h}^{(x-0)/h}
w_{k}^{\alpha} V(x-kh,t) + \mathcal{O}(h^{2})
\end{smallmatrix} \right).
\end{eqnarray}
where
\begin{equation}\label{12a}
w_{k}^{\alpha}=\frac{(-1)^{k}\Gamma(\alpha+1)}{\Gamma(\alpha/2-k+1)\Gamma(\alpha/2+k+1)}.
\end{equation}
Denote by $U_{i}^{h}(t)$  and $V_{i}^{h}(t)$ the approximated  values
 of $U(x_{i},t)$ and $V(x_{i},t)$ respectively,
 the central finite difference approximation is therefore given by
 \begin{eqnarray}
\label{dis11}
\dfrac{\partial^{\alpha}  B^{1}(x_{i},t) }{\partial|x|^{\alpha}}\approx
\left(\begin{smallmatrix}
 -\dfrac{1}{h^{\alpha}} \sum\limits_{k=-M+i}^{i}
w_{k}^{\alpha} U_{i-k}^{h}(t)\\
-\dfrac{1}{h^{\alpha}}\sum\limits_{k=-M+i}^{i}
w_{k}^{\alpha} V_{i-k}^{h}(t)
\end{smallmatrix} \right), i =1,.., M-1.
\end{eqnarray}
Note that  $U_0^h(t)=\text{Re}(\phi_1(t))$, $U_M^h(t)=\text{Re}(\phi_2(t))$, $V_0^h(t)=\text{Im}(\phi_1(t))$ and  $V_M^h(t)=\text{Im}(\phi_2(t))$.
By setting $ U_{h}=(U_{i}^{h})_{1\leq i \leq M-1}$, $ V_{h}=(V_{i}^{h})_{1\leq i \leq M-1}$ and  $B_{h}^{1}=(U_h,V_h)^{T}$, the semi discrete version of \eqref{compac}
after space discretization is given  by
\begin{eqnarray}
 \label{compacd}
 \left\lbrace \begin{array}{l}
 \dfrac{d B_{h}^{1}}{dt}=\mathcal{A}_{h} B_{h}^{1}+ F( B_{h}^{1}),\\
 B_{h}^{1}(0) = \left( g(x_i)\right)_{1\leq i \leq M-1}.
 \end{array}\right.
\end{eqnarray}
 where
 \begin{eqnarray}
\label{new112}
  \mathcal{A}_h= \left(\begin{smallmatrix}  -\mathbf{P}+\gamma_{r} \mathbf{I} & -\gamma_{i} \mathbf{I}\\ \gamma_{i} \mathbf{I} &  -\mathbf{P} +\gamma_{r}\mathbf{I} \end{smallmatrix} \right), \\
  \mathbf{P}=(p_{i\,j})_{1\leq i,j \leq M-1},\,\,\,\, p_{i\,j}= \frac{Pr}{h^{\alpha}} w_{i-j}^{\alpha},\, F( B_{h}^{1})=F_1( B_{h}^{1})+Bc(t),
   \end{eqnarray}
   $Bc(t)$ being the contribution of the Dirichlet boundary condition, which should be expressed as a function of $\phi_1(t)$ and $\phi_2(t)$.
 Note  that $\mathbf{I}$ is  the $(M-1) \times (M-1)$ identity matrix. Please also  note that the matrix $\mathcal{A}_h$ is more than 50 $\%$ full.

 Let $N$ being the time subdivision, we use  the constant time step $\tau=T/N$. In order to  fully discretize \eqref{compac}, let
 $B_{h,n}^{1}$ being our approximated  solution of $B^{1}(n \tau) = \left(B^{1}( x_i, n \tau) \right)_{1 \leq i \leq M-1}$.
From \eqref{compacd}, the $\theta-$ Euler  Riesz fractional finite-difference scheme  to approximate  \eqref{compac} is given by
\begin{eqnarray}
\label{impl}
\dfrac{B_{h,n+1}^{1}- B_{h,n}^{1}}{\tau}=\theta \left( \mathcal{A}_h B_{h,n}^{1} +F( B_{h,n}^{1})\right)+(1-\theta) \left( \mathcal{A}_h B_{h,n+1}^{1}+F( B_{h,n+1}^{1})\right)\\
0\leq\theta\leq 1. \nonumber
\end{eqnarray}
For  $\theta=1$, the scheme is an explicit Riesz fractional
finite-difference scheme, while for  $\theta=0$, the scheme is a
fully implicit Riesz fractional finite-difference scheme.  The
high order accuracy in time is obtained for  $\theta=\frac{1}{2}$,
which corresponds to the Crank-Nicholson Riesz fractional
finite-difference approximation scheme. Note  the  for
$\theta=0$, the corresponding explicit scheme is only stable for
very small time step $\tau$. For $ \theta \neq 1$, the scheme is
more stable, but the fact that the matrix $\mathcal{A}_h$   is likely to be  more than  50 $\% $ full (depending of the sparseness of  $\mathbf{P}$) makes the Newton iterations less efficient. High order implicit Runga-kutta methods can be used if
high order accuracy is needed, but  these methods will  be extremely less efficient.

To solve the efficiency drawback of the  implicit schemes, we
propose in this work two simple schemes.

For  nonstiff nonlinear  part $F$,   we consider the semi-implicit Riesz fractional
finite-difference scheme where the linear part of \eqref{compacd}
is approximated implicitly and the nonlinear part explicitly.
Following \cite{ATthesis,Antoine,TLGspe}, the corresponding scheme is
given by
\begin{eqnarray}
\label{semi}
\dfrac{B_{h,n+1}^{1}- B_{h,n}^{1}}{\tau}=  \mathcal{A}_h B_{h,n+1}^{1} +F( B_{h,n}^{1}).
\end{eqnarray}
For  stiff nonlinear  part $F$, scheme \eqref{semi} will require small time steps to be stable, following \cite{Antoingajan}
 we consider the following  the linearly  implicit Riesz fractional
finite-difference scheme given  by
\begin{eqnarray}
\label{lemi}
B_{h,n+1}^{1}&=& B_{h,n}^{1} + \tau \left(\mathbf{I}-\tau J_n\right)^{-1}\left(\mathcal{A}_h B_{h,n}^{1} +F( B_{h,n}^{1}\right),\\
J_n&=& \mathcal{A}_h+\partial_{B^1} F(B_{h,n}^{1}). \nonumber\\
\end{eqnarray}
Obviously the semi-implicit  scheme  given at \eqref{semi}  and  the linearly  implicit Riesz fractional
finite-difference scheme are very efficient than the  implicit schemes given in \eqref{impl} for $ \theta \neq 1$,  as only one linear system is solved per time iteration.

Following again \cite{Antoingajan}, we can also  obtain  the $s-$stages Rosenbrock Riesz fractional
finite-difference schemes  if high order accuracy is  needed.  Such schemes will  be  efficient  as only $s$ linear systems are required by time iteration.

\subsection{Numerical results}
\subsubsection{Numerical simulations of  the   CFGL equation  and localized wave solutions}
Numerical  simulations of \eqref{15} are performed using the
the linearly implicit  Riesz fractional finite-difference scheme given by
\eqref{lemi} as  the nonlinear function $F$ is stiff.
We choose the solution of $B^1(x,0)$ in the form of a
nonlinear solution of the standard complex Gingburg Landau equation
\cite{Nozaki,Stenflo}

\begin{equation}\label{21}
B^{1}(x,0)= \frac{B_{0}e^{\theta}}{1+e^{(\theta+\theta^{*})^{(1+i\mu)}}}, \,\,\phi_1(t)=\phi_2(t)=0,\,\, t \in (0.T].
\end{equation}
where the real part and the imaginary part of $B^{1}(x,0)$ are
given respectively by

\begin{eqnarray}
\label{35} \left\lbrace \begin{array}{l}
B_{r}^{1}(x,0)=U(x,0)=B_{0}\Big[\frac{e^{-\theta}+\cos(2\mu
\theta)e^{\theta}}{2(\cosh(2\theta)+\cos(2\mu \theta))}\Big],\\
B_{i}^{1}(x,0)=V(x,0)= -B_{0}\Big[\frac{\sin(2\mu
\theta)e^{\theta}}{2(\cosh(2\theta)+\cos(2\mu \theta))}\Big],
\end{array}\right.
\end{eqnarray}
where $\theta=kx$, $\mu=-\beta+\sqrt{2+\beta^2}$ and
$\beta=-\frac{3Q_{r}}{2Q_{i}}$.

The parameter values are: $\Omega_{0}^{2}=0.032$, $k=1.5$, $I=3$
$\lambda_{1}=2$, $\lambda_{3}=0.992$, $\eta_{0}=1$, $\eta_{1}=6$,
$\eta_{2}=3$, $r=0.008$, $c_{0}=0.001$, $c_{1}=0.001$  and $B_{0}=
0.5$.

\begin{figure}[h!]
\subfigure[]{
\includegraphics[width=4.3in]{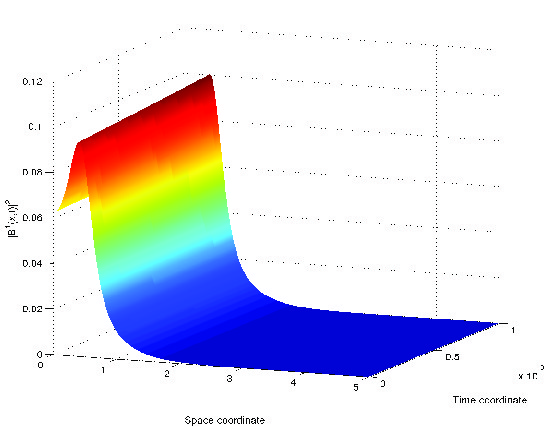}}
\caption{(Color online) Profile of the numerical solution
$|B^1(x,t)|^2$ according to time and space for $\alpha=1.8$  at up
to the final time  $T=0.001$. }\label{fig2}
\end{figure}

\figref{fig2} displays the spatiotemporal evolution on the amplitude
for $\alpha=1.8$ at time $t=0.001$.
We observe in this figure that the solution is well localized nonlinear excitation in space and time,
it  has the shape of a short pulse and propagates without any change
of its profile. It is clear from there that as time evolves, the
form of the pulse does not change; it is structurally stable.

\begin{figure}[h!]
\centering
\includegraphics[width=4.3in]{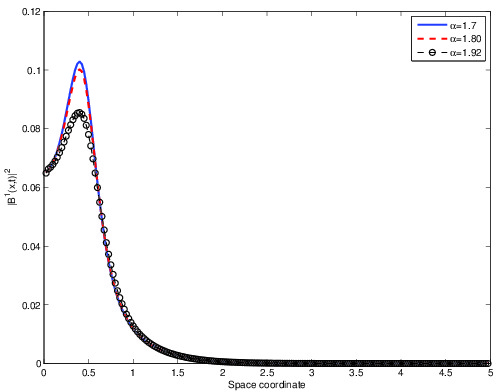}
\caption{Spatial profile of the numerical solution $|B^1(x,t)|^2$ at
time $t=0.001$ with initial nonlinear localized solution with
$\alpha=1.7$,  $\alpha=1.8$, and  $\alpha=1.92$.} \label{fig3}
\end{figure}

\figref{fig3} displays spatial profiles of the amplitude of the
solution for three distinct values of parameter $\alpha$, namely
$\alpha=1.7$, $\alpha=1.8$ and  $\alpha=1.92$  for the time instant
$t=0.001$. We observe in the graphs of \figref{fig3} that the
solution is well localized in space with the shape of a short pulse
and it amplitude decreases with the increasing of $\alpha$.
 Then, The fractional order mostly contributes to the
behavior of the tails of the short pulse numerical solutions.
Remarkably, the pulse profiles in \figref{fig3} are in qualitative
agreement with the typical results reported in electrodynamics
theory in both myelinated and myelin-free nerve fiber contexts
\cite{Dikande}.

\subsubsection{Reconstruction of  the discrete solutions from continuous CFGL solution}

In order to check the validity of our fractional approach and to get
an idea of what kind of dynamical waves one might obtain in the
neural network, we carried out numerical simulation of \eqref{1} and
compare  the results with the one of the fractional model
\eqref{new}. The simulation  of \eqref{1} is performed through the
fourth-order Runge-Kutta scheme.

Remember that  $$B^{1}(x,t)=B_{r}^{(1)} (x,t)+i B_{i}^{(1)}(x,t)=U(x,t)+i V(x,t).$$
From the first equation of \eqref{5}, we have
\begin{equation}\label{36}
\psi=2\varepsilon(B_{r}^{(1)}\cos\theta-B_{i}^{(1)}\sin\theta)
+\varepsilon^{2}[C^{(1)}+2(D_{r}^{(1)}\cos2\theta-D_{i}^{(1)}\sin
2\theta)],
\end{equation}
where  $\theta =k\,x - \Omega t$, $D_{r}^{(1)}$
 and $D_{i}^{(1)}$  are the real and imaginary parts of
$D^{(1)}$.
From  \eqref{11}  we have
\begin{equation}\label{37}
D^{(1)}=(a_{1}-ia_{2})[B^{1}]^{2},\,\qquad
a_{1}= \frac{\lambda_{1}}{\Omega_{1}^{2}}, \,\qquad a_{2}=
\frac{\Omega\eta_{1}}{\Omega_{1}^{2}},
\end{equation}
which leads to
\begin{eqnarray}\label{12a}
D_{r}^{(1)}(x,t)&=&a_{1}(B_{r}^{(1)2}-B_{i}^{(1)2})+2a_{2}B_{r}^{(1)}B_{i}^{(1)},\\
D_{i}^{(1)}(x,t)&=&a_{2}(B_{i}^{(1)2}-B_{r}^{(1)2})+2a_{1}B_{r}^{(1)}B_{i}^{(1)}.
\end{eqnarray}

Inserting \eqref{37} into \eqref{36} and using the relation
$u=\varepsilon\psi$, we obtain for the nerve impulse the following
solution
\begin{eqnarray}\label{42}
\begin{split}
u &=2\varepsilon^{2}\Big\{[B_{r}^{(1)}\cos(\theta )-B_{i}^{(1)}\sin(\theta)]
-\frac{\lambda_{1}}{\Omega_{0}^{2}}[B_{r}^{(1)2}+B_{i}^{(1)2}]\Big\}\\
&+2\varepsilon^{3}\Big\{[a_{1}(B_{r}^{(1)2}-B_{i}^{(1)2})
+2a_{2}B_{r}^{(1)}B_{i}^{(1)}]\cos(2\theta )\\
&+[a_{2}(B_{i}^{(1)2}-B_{r}^{(1)2})
+2a_{1}B_{r}^{(1)}B_{i}^{(1)}]\sin(2\theta )\Big\}.
\end{split}
\end{eqnarray}
Applying a similar procedure, we obtain for the bursting variable
the following initial solution
\begin{equation}\label{43}
\Phi=2\varepsilon(B_{r}^{(2)}\cos\theta-B_{i}^{(2)}\sin\theta)
+\varepsilon^{2}[C^{(2)}+2(D_{r}^{(2)}\cos2\theta-D_{i}^{(2)}\sin
2\theta)],
\end{equation}
where $B_{r}^{(2)}$ ($D_{r}^{(2)}$) and $B_{i}^{(2)}$
($D_{r}^{(2)}$) are the real and imaginary parts of $B^{(2)}$
($D^{(2)}$), respectively. Then we have
\begin{equation}\label{44}
\begin{split}
v &=\varepsilon \Phi\\
&=2\varepsilon^{2}(B_{r}^{(2)}\cos\theta-B_{i}^{(2)}\sin\theta)
+\varepsilon^{3}[C^{(2)}+2(D_{r}^{(2)}\cos2\theta-D_{i}^{(2)}\sin
2\theta)],
\end{split}
\end{equation}
where

\begin{equation}\label{45}
B_{r}^{(2)}=B_{i}^{(2)}=0,
\end{equation}
\begin{equation}\label{46}
C^{(2)}=-\frac{2d}{e}[B_{r}^{(1)2}+B_{i}^{(1)2}],
\end{equation}
\begin{equation}\label{47}
D_{r}^{(2)}=-\frac{d}{e}[B_{r}^{(1)2}-B_{i}^{(1)2}],
\end{equation}
\begin{equation}\label{48}
 D_{i}^{(2)}=-\frac{2d}{e}[B_{r}^{(1)}\times B_{i}^{(1)}].
\end{equation}

We also have
\begin{equation}\label{49}
\beta=2\varepsilon(B_{r}^{(3)}\cos\theta-B_{i}^{(3)}\sin\theta)
+\varepsilon^{2}[C^{(3)}+2(D_{r}^{(3)}\cos2\theta-D_{i}^{(3)}\sin
2\theta)],
\end{equation}
where $B_{r}^{(3)}$ ($D_{r}^{(3)}$) and $B_{i}^{(3)}$
($D_{i}^{(3)}$) are the real and imaginary parts of $B^{(3)}$
($D^{(3)}$), respectively. Then we have
\begin{equation}\label{50}
\begin{split}
w &=\varepsilon \beta\\
&=2\varepsilon^{2}(B_{r}^{(3)}\cos\theta-B_{i}^{(3)}\sin\theta)
+\varepsilon^{3}[C^{(3)}+2(D_{r}^{(3)}\cos2\theta-D_{i}^{(3)}\sin
2\theta)],
\end{split}
\end{equation}
where
\begin{equation}\label{51}
B_{r}^{(3)}(x,t)=b_{1}B_{r}^{(1)}-b_{2}B_{i}^{(1)},
\end{equation}
\begin{equation}\label{52}
B_{i}^{(3)}(x,t)=b_{1}B_{i}^{(1)}+b_{2}B_{r}^{(1)},
\end{equation}
\begin{equation}\label{53}
D_{r}^{(3)}(x,t)=c_{11}D_{r}^{(1)}-c_{2}D_{i}^{(1)},
\end{equation}
\begin{equation}\label{54}
D_{i}^{(3)}(x,t)=c_{11}D_{i}^{(1)}+c_{2}D_{r}^{(1)},
\end{equation}

\begin{equation}\label{55}
C^{(2)}(x,t)=\frac{\Omega_{0}^{2}}{r}C^{(1)},
\end{equation}
with
\begin{equation}\label{56}
\begin{split}
&b_{1}= \frac{\Omega_{0}^{2}r}{r^{2}+\Omega^{2}}, \;\
\;\ b_{2}= \frac{\Omega_{0}^{2}\Omega}{r^{2}+\Omega^{2}},\\
&c_{11}= \frac{\Omega_{0}^{2}r}{r^{2}+4\Omega^{2}}, \;\ \mathrm{and}
\;\ c_{2}= \frac{2\Omega_{0}^{2}\Omega}{r^{2}+4\Omega^{2}}.
\end{split}
\end{equation}

To have discrete  solutions from  continuous fractional model \eqref{new}, we proceed as follows
\begin{itemize}
 \item We use the initial solution \eqref{35} at  the  points $ \theta=\theta_{n,0}=k \,x_n$, where $x_n= n h,\, 1\leq n\leq M-1$, and obtain from the linearly implicit Riesz fractional
finite-difference  scheme the approximated solution $B_{h,m}^{1}\backsimeq \left( B^{1}(x_n, t_m)\right)_{1\leq n\leq M-1}$.
\item The numerical solution $B_{h,m}^{1}$ is now used in \eqref{42},\eqref{44} and \eqref{50} with $\theta=\theta_{n,m}=k x_n-\Omega t_m$ to obtain  the discrete solution
$u_n(t_m):=u(x_{n},t_m)$,$v_n(t_m):=v(x_{n},t_m)$ and $w_n(t_m):=w(x_{n},t_m)$.
\end{itemize}
The  discrete solutions $u_n(t_m)$,$v_n(t_m)$ and $w_n(t_m)$ obtained from the continuous fractional model \eqref{new}, can be compared with the numerical solution of \eqref{1} from implicit fourth-order
Runge-Kutta method. Note that the initial solution used to solve \eqref{1} is obtained from \eqref{42},\eqref{44} and \eqref{50}  with the same initial solution \eqref{35} for
$\theta=\theta_{n,0}=k x_n$.

In our graphs,  the surface plots of $u_n$,$v_n $ and $w_n$ from the
continuous fractional model \eqref{new} will be called fractional
discrete solutions, while the one coming directly from \eqref{1}
will  be called discrete solutions. We should always remember that
to have the fractional model \eqref{new} the time $t$, the
coefficients $\eta_0$, $c_1$ and $\lambda_3$ have be perturbed at
order $\varepsilon^2$. So, the fractional discrete  solution at time
$t$ will be compared with the discrete solution at time
$\dfrac{t}{\varepsilon^2}$.

\begin{figure}[!ht]
\subfigure[]{
 \includegraphics[width=3in]{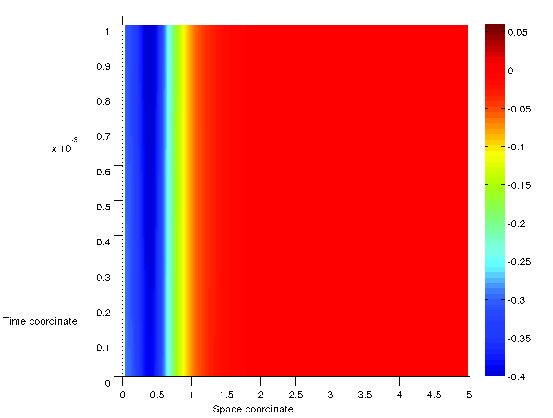}}
\hskip 0.01\textwidth
\subfigure[]{
 \includegraphics[width=3in]{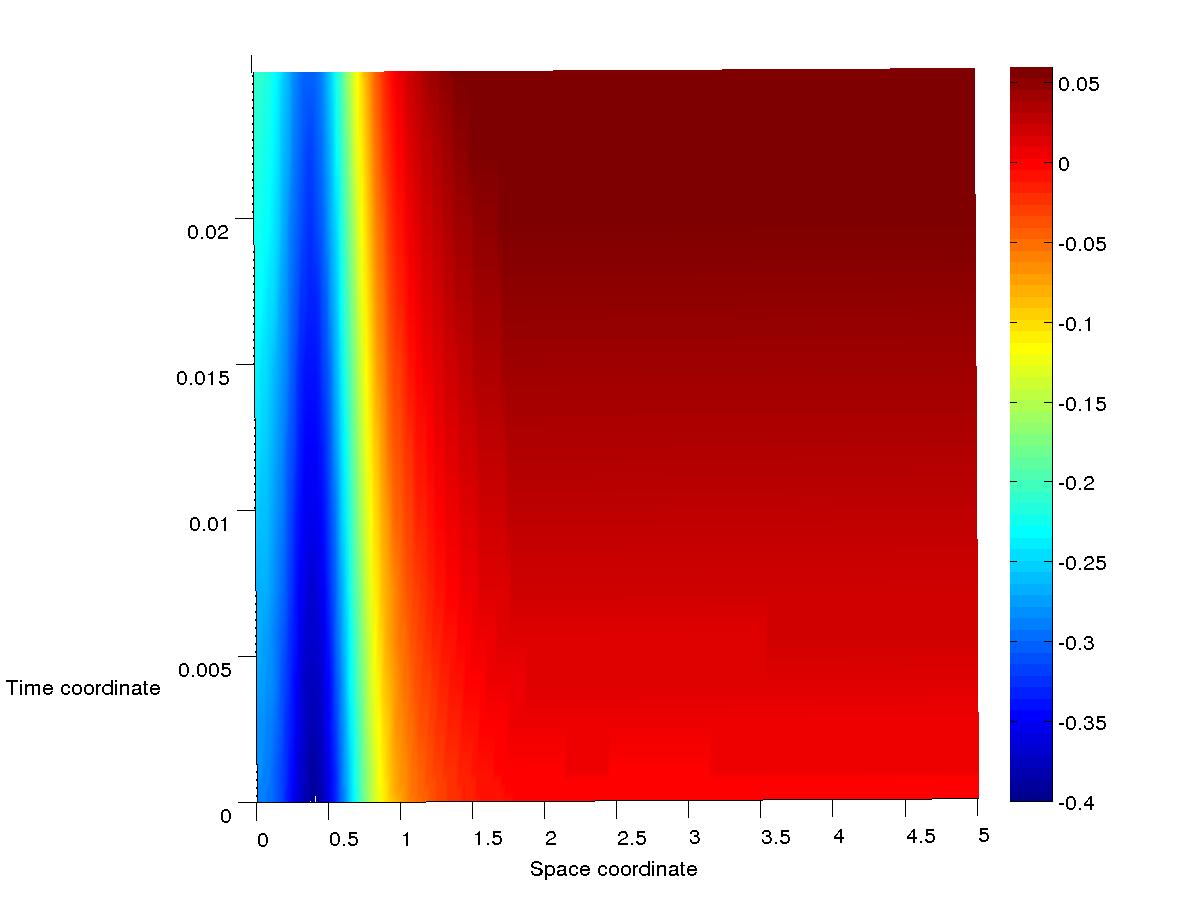}}
 \subfigure[]{
 \includegraphics[width=3in]{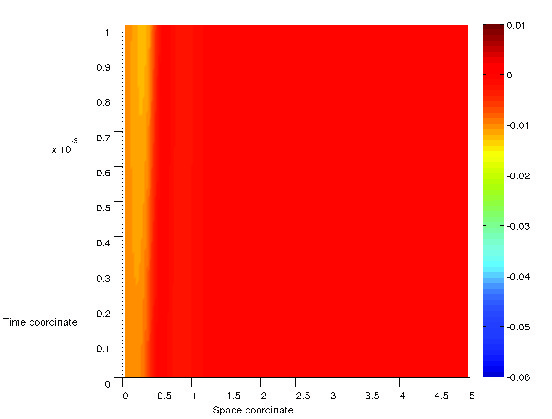} }
\hskip 0.01\textwidth
\subfigure[]{
 \includegraphics[width=3in]{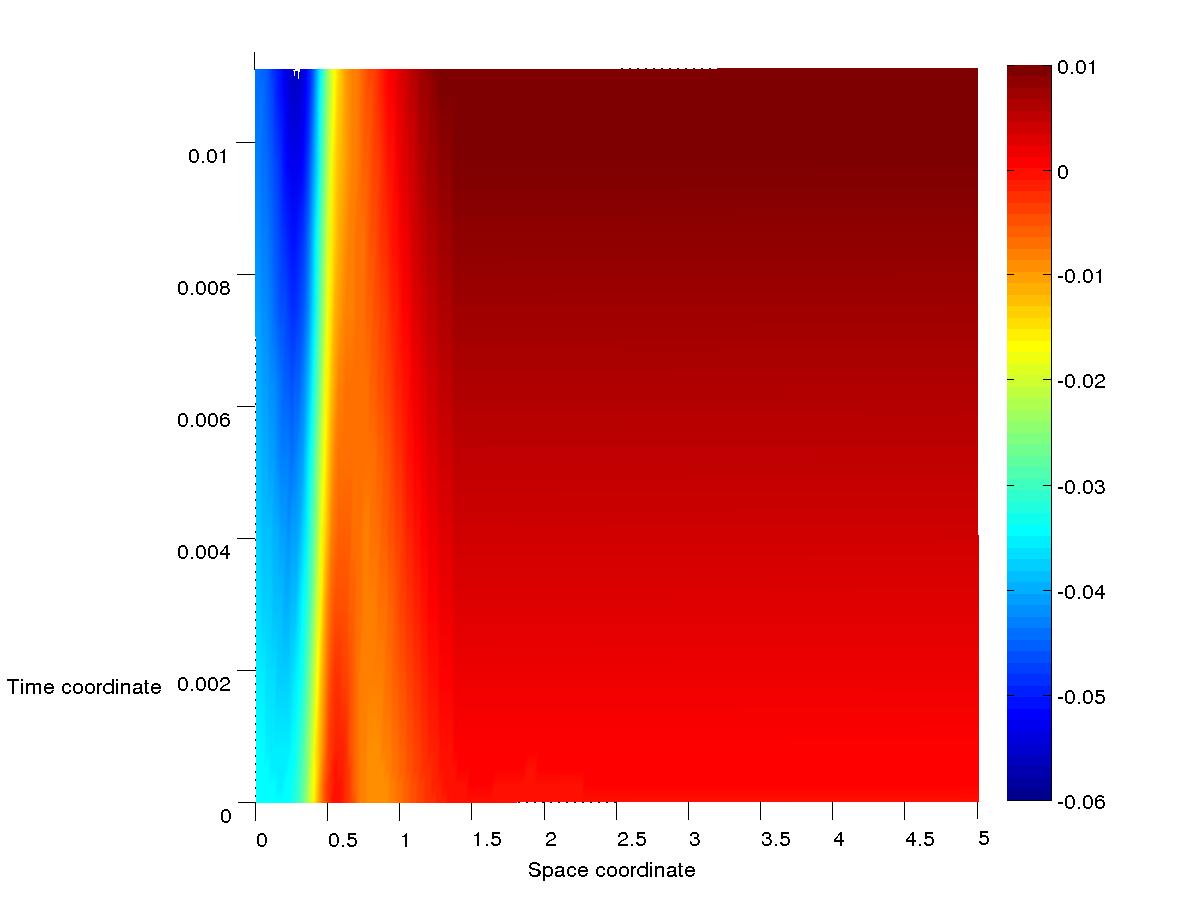}}
 \subfigure[]{
 \includegraphics[width=3in]{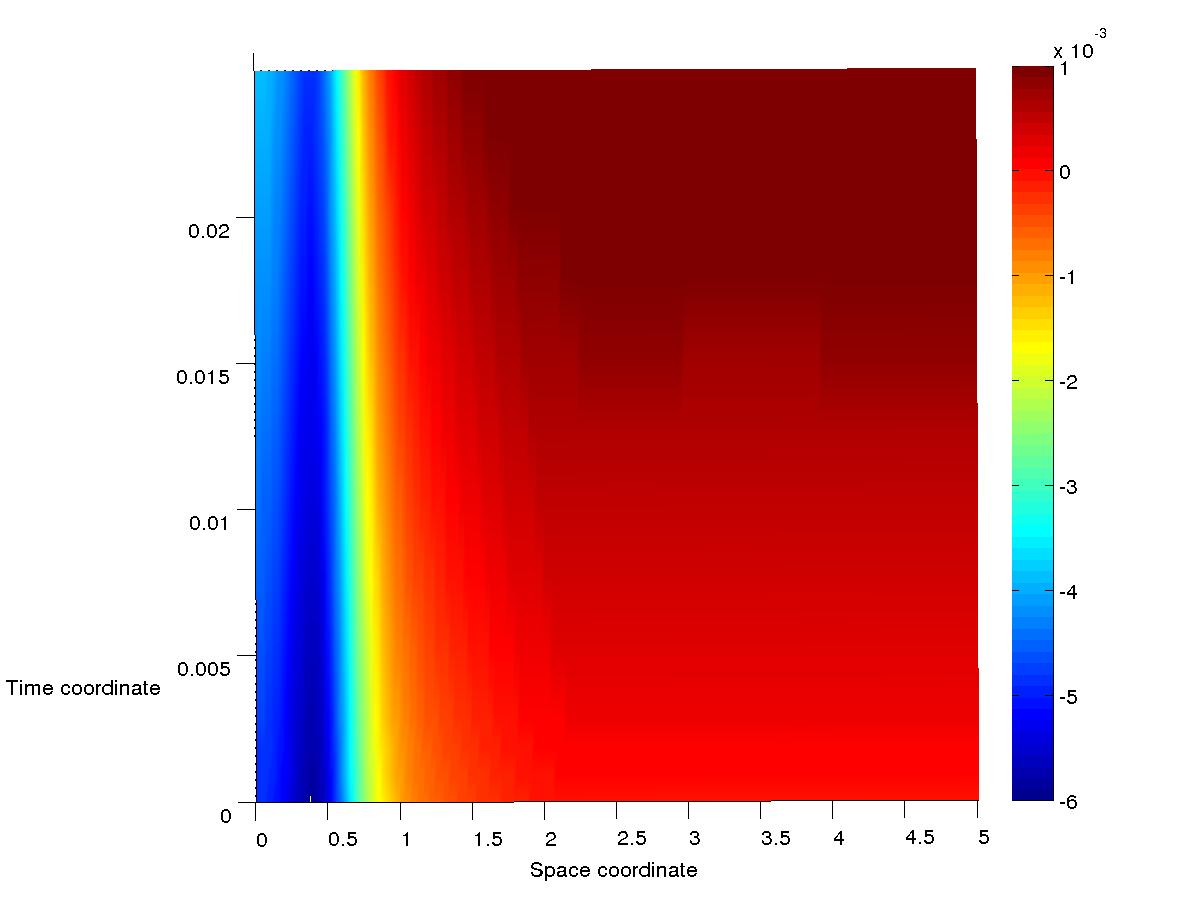}}
\hskip 0.01\textwidth
\subfigure[]{
 \includegraphics[width=3in]{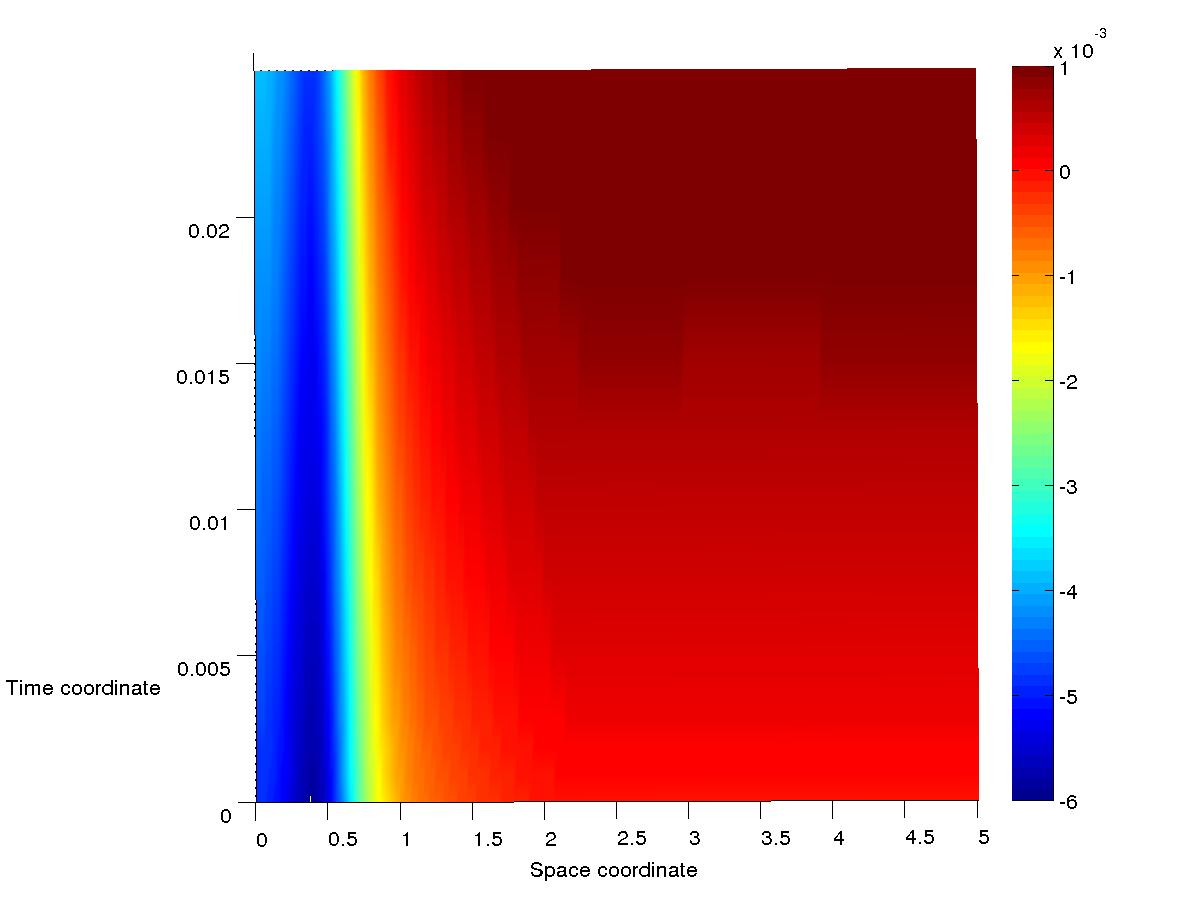}}
\caption{ Fractional $u_n$ (a) and  discrete $u_n$(b),  fractional $v_n$ (c) and discrete $v_n$ (d), and  fractional $w_n$ (e) and discrete $w_n$ (f).
The fractional solutions are up to final time $T=0.001$  while  the discrete solutions are up to the final time $ T_f=\frac{T}{\varepsilon^2}=0.025$ with $\varepsilon=0.2$}
\label{FIG04}
\end{figure}

\begin{figure}[!ht]
\subfigure[]{
 \includegraphics[width=3in]{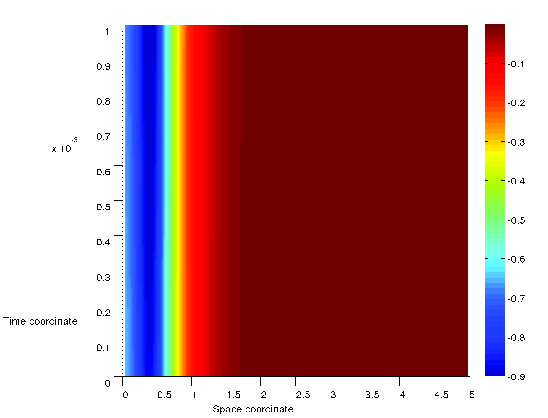}}
\hskip 0.01\textwidth
\subfigure[]{
\includegraphics[width=3in]{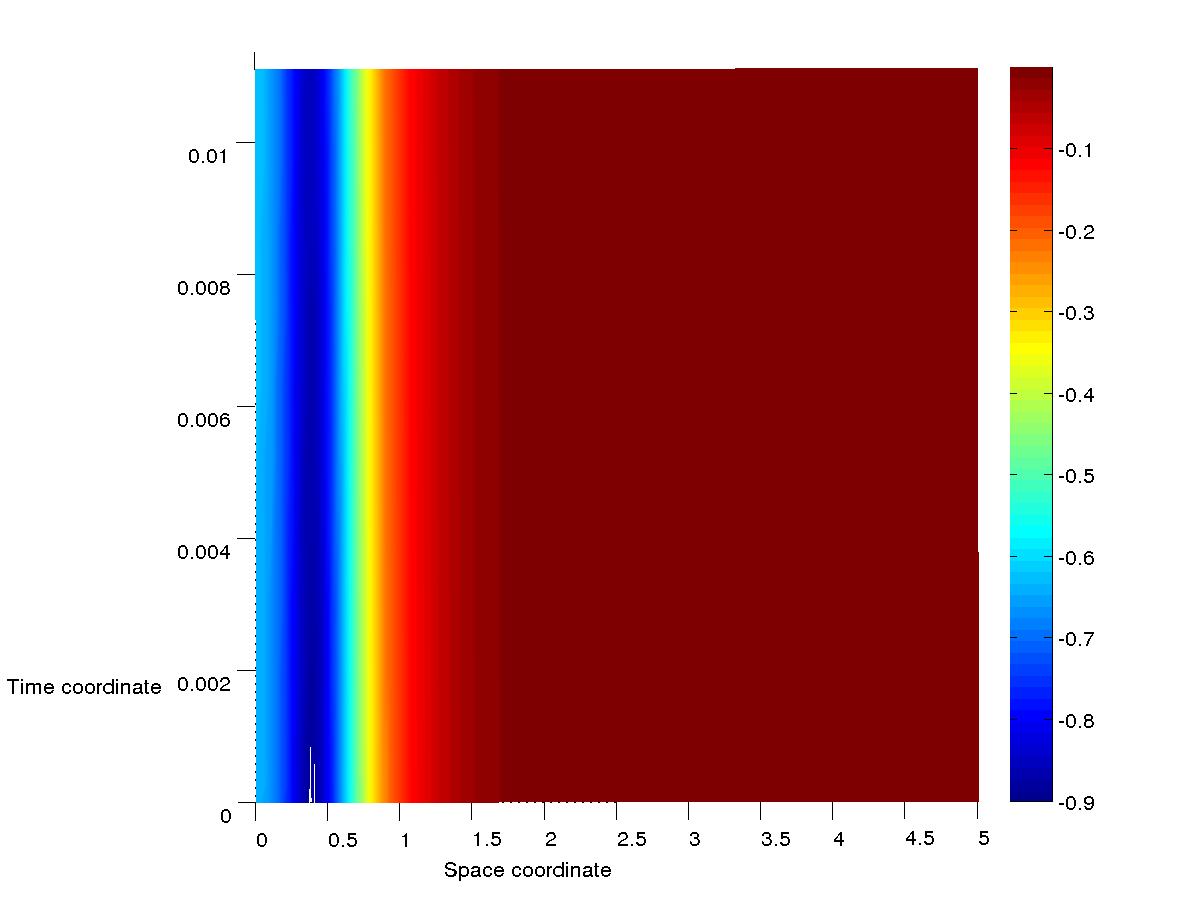}}
 \subfigure[]{
 \includegraphics[width=3in]{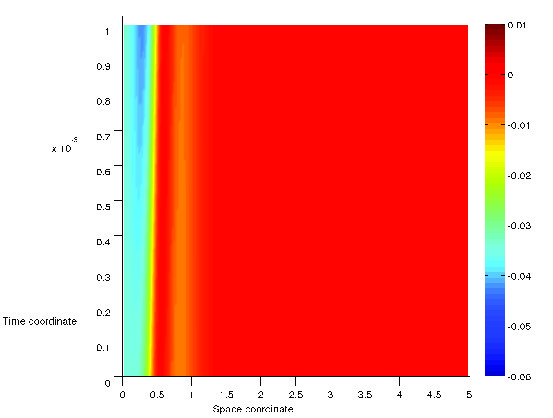} }
\hskip 0.01\textwidth
\subfigure[]{
 \includegraphics[width=3in]{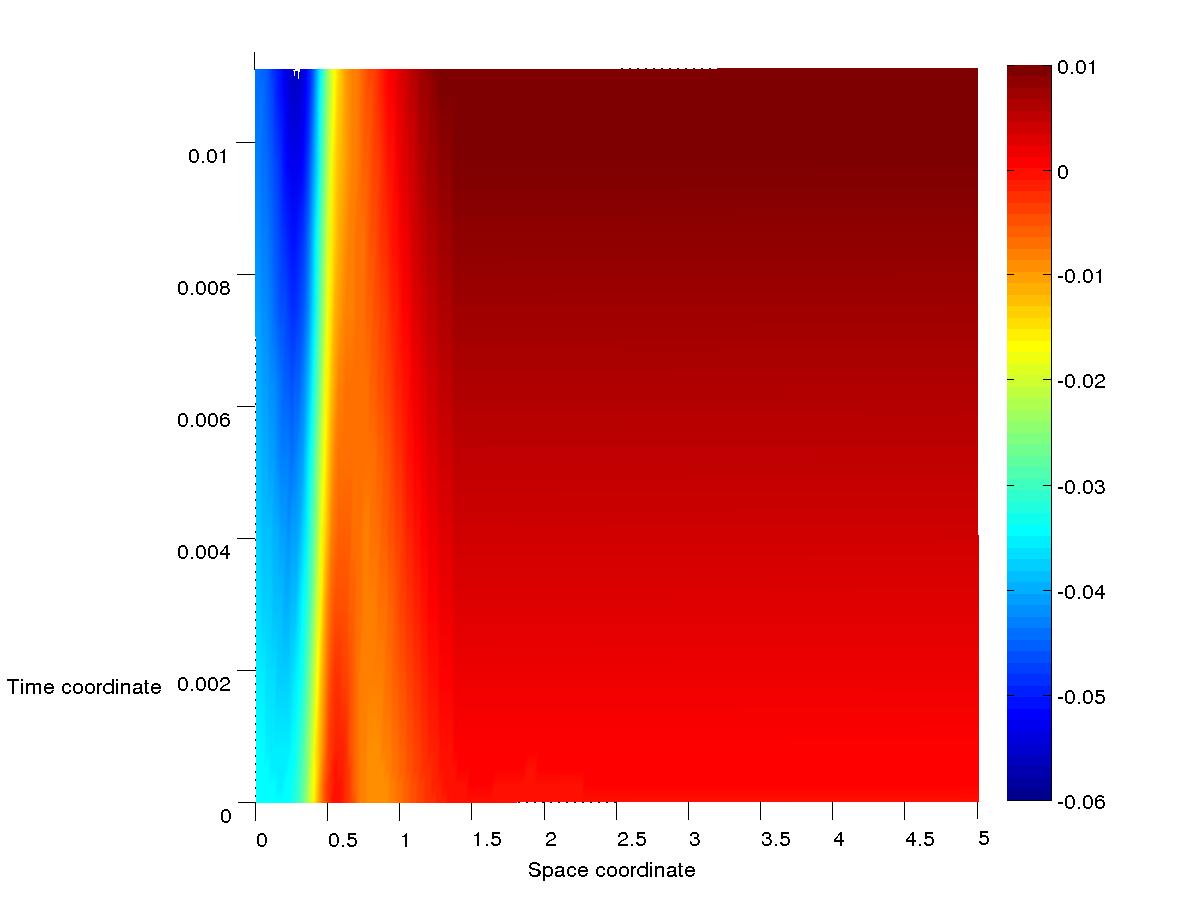}}
 \subfigure[]{
 \includegraphics[width=3in]{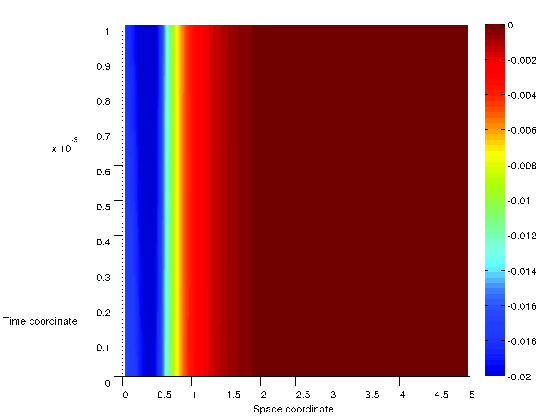}}
\hskip 0.01\textwidth
\subfigure[]{
 \includegraphics[width=3in]{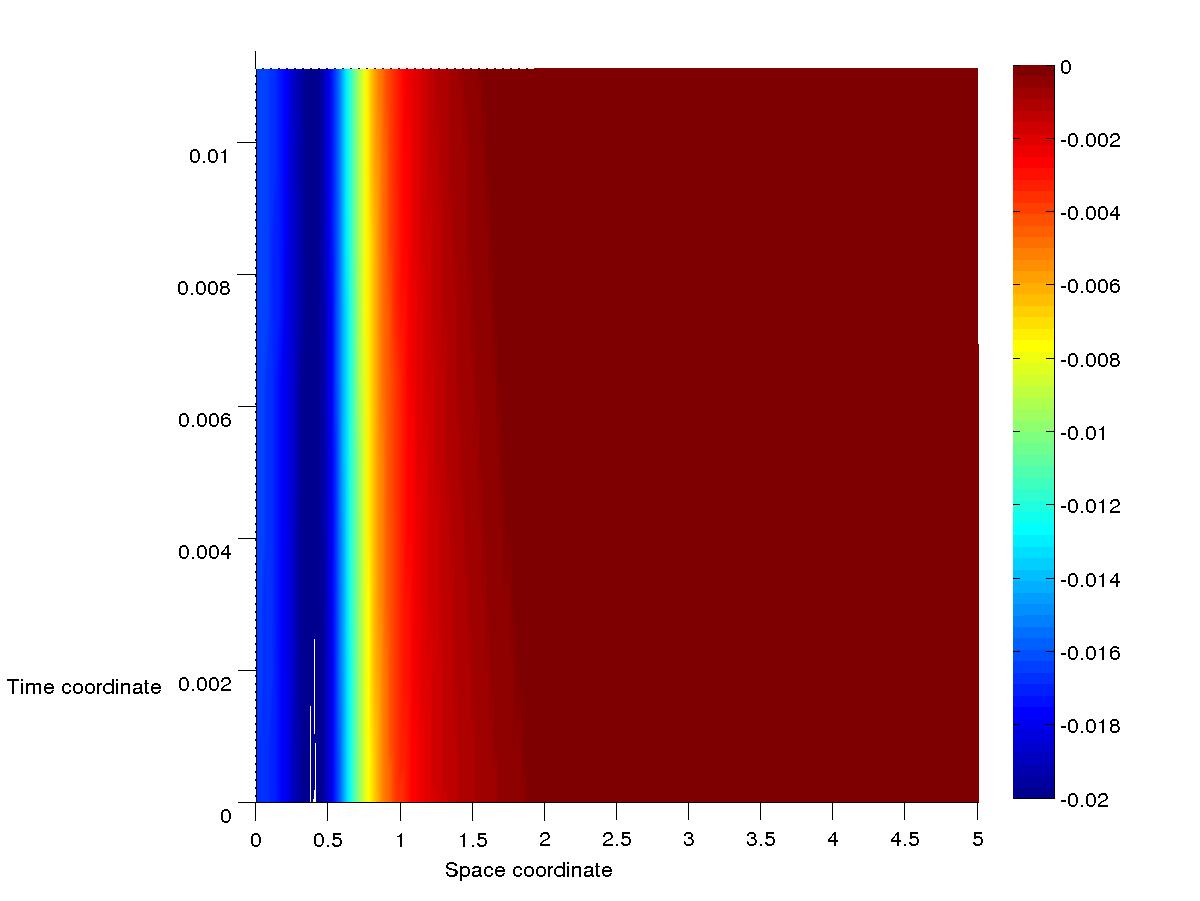}}
\caption{ Fractional $u_n$ (a) and  discrete $u_n$ (b),  fractional $v_n$ (c) and discrete $v_n$ (d), and  Fractional $w_n$ (e) and discrete $w_n$ (f).
The fractional solutions are up to final time $T=0.001$  while  the disrete solutions are up to the final time $ T_f=\frac{T}{\varepsilon^2}=0.0111$ with $\varepsilon=0.3$.}
\label{FIG05}
\end{figure}

\begin{figure}[!ht]
\subfigure[]{
 \includegraphics[width=3in]{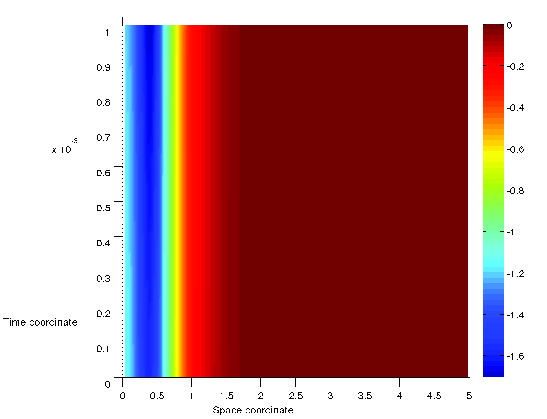}}
\hskip 0.01\textwidth
\subfigure[]{
 \includegraphics[width=3in]{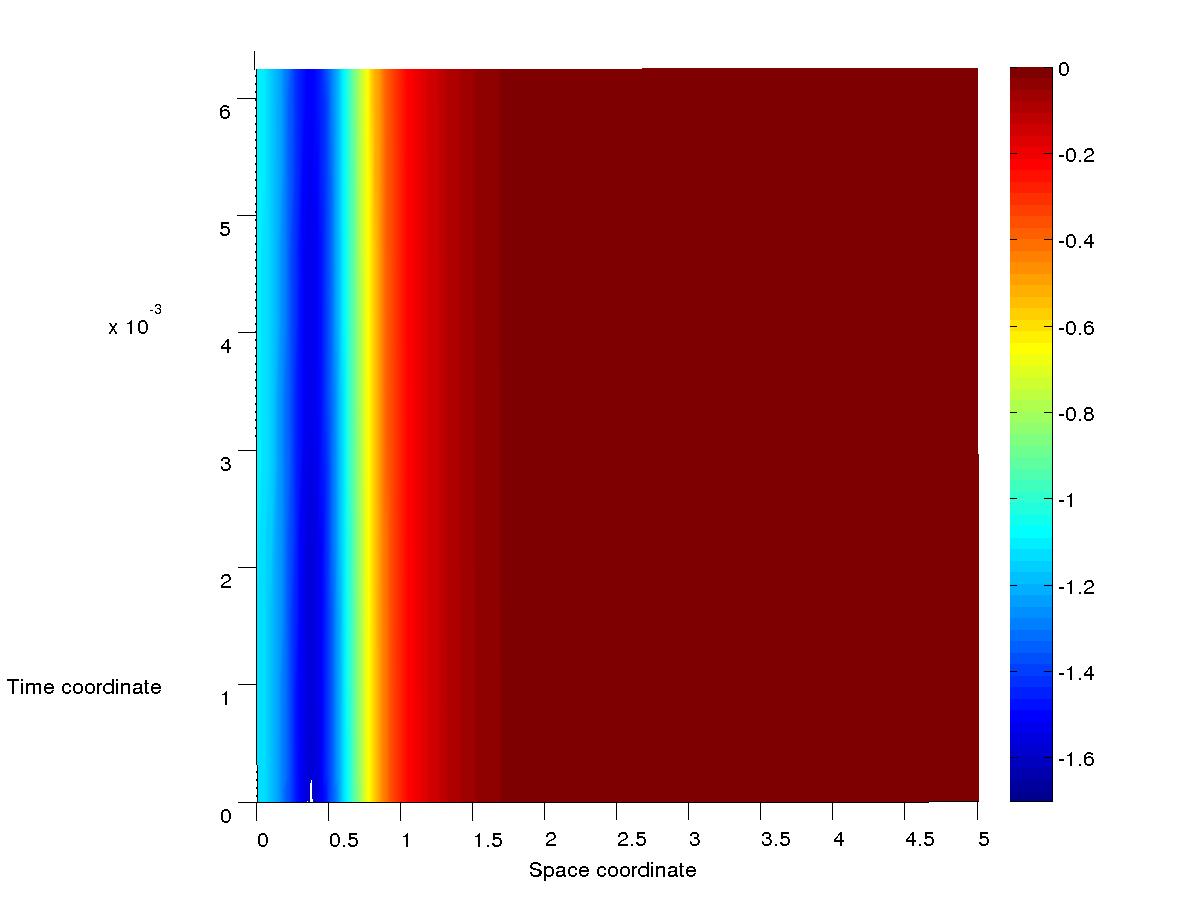}}
 \subfigure[]{
 \includegraphics[width=3in]{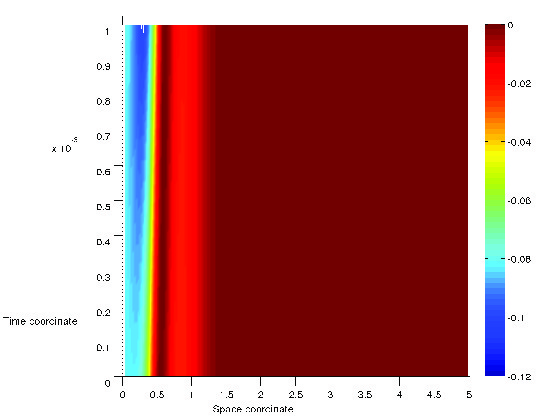} }
\hskip 0.01\textwidth
\subfigure[]{
 \includegraphics[width=3in]{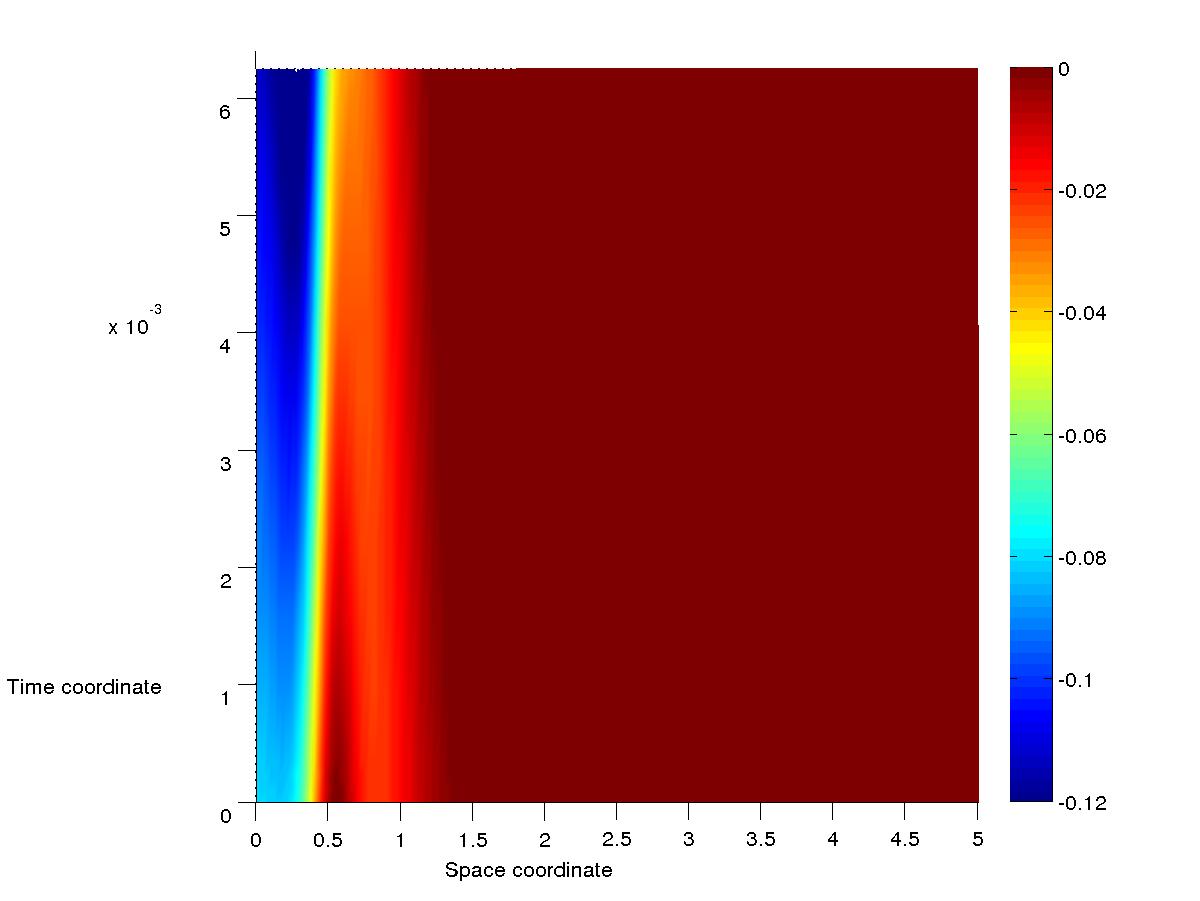}}
 \subfigure[]{
 \includegraphics[width=3in]{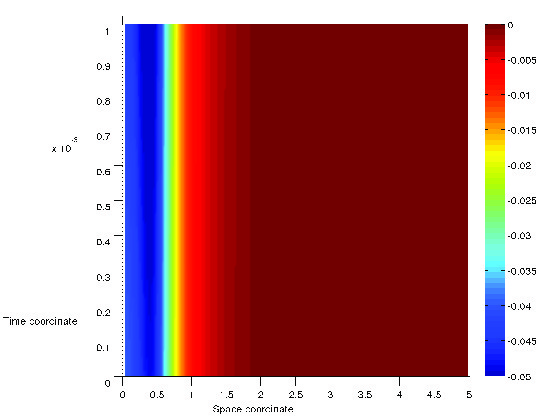}}
\hskip 0.01\textwidth
\subfigure[]{
 \includegraphics[width=3in]{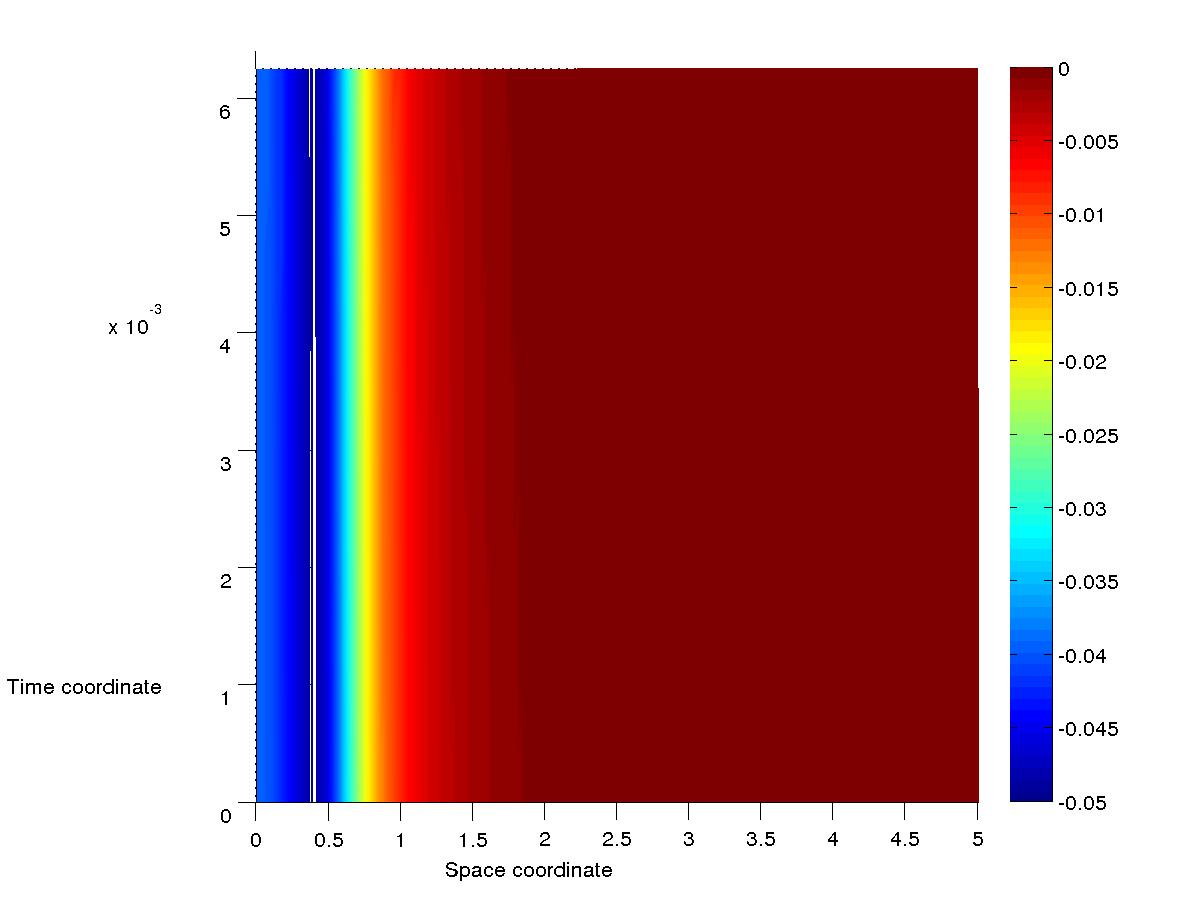}}
\caption{ Fractional $u_n$ (a) and  discrete $u_n$(b),  fractional $v_n$ (c) and discrete $v_n$ (d), and  fractional $w_n$ (e) and discrete $w_n$ (f).
The fractional solutions are up to final time $T=0.001$  while  the discrete solutions are up to the final time $ T_f=\frac{T}{\varepsilon^2}=0.0062$ with $\varepsilon=0.4$.}
\label{FIG06}
\end{figure}

\figref{FIG04},  \figref{FIG05} and \figref{FIG06} show the surface
plots of $u_n$,$v_n $ and $w_n$ with different values of
$\varepsilon$. The fractional graphs of  $u_n$  are those  in (a)
and  the discrete graphs of  $u_n$  are those in (b).  The
fractional graphs of  $v_n$  are those in (c) and the discrete
graphs of $v_n$ in (d). The fractional graphs of $w_n$  are those
in (e) and discrete graphs  of $w_n$ in (f).
 All fractional numerical solutions are up to final time $T=0.001$, but in \figref{FIG04}
the discrete numerical  solutions are up to the final time $
T_f=\frac{T}{\varepsilon^2}=0.0111$ with $\varepsilon=0.2$. In
\figref{FIG05}, the discrete numerical  solutions are up to the
final time $ T_f=\frac{T}{\varepsilon^2}=0.0111$ with
$\varepsilon=0.3$, while in \figref{FIG06}, the  numerical discrete
solutions are up to final  time $
T_f=\frac{T}{\varepsilon^2}=0.0062$ with $\varepsilon=0.4$. From
those graphs, we can observe that as the final time  $
T_f=\frac{T}{\varepsilon^2}$ decreases,  the fractional numerical
solutions of $u_n$,$v_n$ and $w_n$ are extremely close to the
discrete numerical solutions,  therefore the fractional continuous
model at \eqref{new} and  discrete model at \eqref{1} are equivalent
for relatively small time.  The long-range effect therefore mostly
contributes to the behavior of the tails of the solutions, as
already mentioned.

\section{Conclusion}
\label{sec5}

The goal of this paper was to study  the nonlinear dynamics of a
diffusively Hindmarsh-Rose neural network with long-range couplings.
Performing a perturbation technique, we have shown that the dynamics
of modulated waves in our neural network can be elegantly described
by the complex fractional Ginzburg-Landau equation. In general,
exact analytical solutions of fractional nonlinear equations cannot
be obtained. We have proposed according to the stiffness of the
system, the semi implicit Riesz fractional finite-difference scheme
and the linearly implicit Riesz fractional finite-difference scheme
to solve efficiently the complex fractional Ginzburg-Landau
equation. It has been revealed that the numerical solutions for the
nerve impulse are well-localized stable short impulses. The results
have been confirmed by the numerical simulations of discrete
equations. The work suggests that long-range diffusive couplings
could be perceived as a way to transport information via nonlinear
waves both in spatial and temporal dimensions for specific processes
of the brain to be controlled. The fractional properties observed in
our neural network may be also advantageous in excitable systems for
crucial intuitions into spatio-temporal dynamics, synchronization
and chaos. The work gives also the opportunity to familiarize with
improved fractional analytical and numerical methods which can be
used to study other systems with long-range couplings.


\section*{Acknowledgments}

A. Mvogo thanks Dr. habil. Anatole Kenfack  of the Institute of
Chemistry and Biochemistry, Freie Universitaet Berlin for fruitful
discussions. A. Tambue was supported by the Robert Bosch Stiftung
within the ARETE chair programme.


\newpage


\end{document}